\documentclass{amsart}
\usepackage{amsmath,amssymb,amsfonts}
\usepackage[mathscr]{eucal}

\usepackage{enumerate}

\usepackage{tikz}
\usepackage{tikz-3dplot}
\usepackage{tikz-cd}
\usetikzlibrary{decorations.markings,decorations.pathmorphing,shapes}
\usepackage{xcolor}
\usepackage{lipsum}

\theoremstyle{plain}
\newtheorem{theorem}{Theorem}[section]

\newtheorem{prop}[theorem]{Proposition}

\theoremstyle{definition}
\newtheorem{definition}[theorem]{Definition}
\newtheorem{example}[theorem]{Example}

\newtheorem*{thank}{Acknowledgments}

\numberwithin{equation}{section}

\input xy
\xyoption{all}

\newcommand{\Ab}{\mathcal Ab}
\newcommand{\Cat}{\mathcal Cat}
\newcommand{\Deltaop}{\Delta^{\op}}
\newcommand{\Gp}{\mathcal Gp}
\newcommand{\Hom}{\operatorname{Hom}}
\newcommand{\id}{\operatorname{id}}
\newcommand{\Map}{\operatorname{Map}}
\newcommand{\nerve}{\operatorname{nerve}}

\newcommand{\op}{\operatorname{op}}
\newcommand{\Sets}{\mathcal Set}
\newcommand{\SSets}{\mathcal{SS}et}
\newcommand{\Top}{\mathcal Top}

\begin{document}

\title[Simplicial sets]{Simplicial sets in topology, category theory, and beyond}

\author[J.E. Bergner]{Julia E. Bergner}

\address{Department of Mathematics, University of Virginia, Charlottesville, VA 22904}

\email{jeb2md@virginia.edu}

\date{\today}

\subjclass[2000]{}


\thanks{The author was partially supported by NSF grant DMS-1906281.}

\begin{abstract}
The notion of a simplicial set originated in algebraic topology, and has also been utilized extensively in category theory, but until relatively recently was not used outside of those fields.  However, with the increasing prominence of higher categorical methods in a wide range of applications, it is important for researchers in a range of fields to have a good working knowledge of them.  This paper is intended as an introduction to simplicial sets, both as an overview of their development from other concepts, and as a user's guide for someone wanting to read modern literature that makes use of them.   
\end{abstract}

\maketitle

\section{Introduction}

Simplicial sets have long been an important structure in the fields of algebraic topology and category theory.  They first developed as a means to model topological spaces purely combinatorially, and as such have been used extensively in homotopy theory.  For many topologically-minded users of simplicial sets, the relationship between the two structures is so entwined that simplicial sets are simply referred to as ``spaces".  On the other hand, simplicial sets are most elegantly described in category-theoretic terms, and via the nerve construction small categories can themselves be described in terms of simplicial sets.  

More recently, the utility of simplicial sets has expanded beyond these fields to a wide range of applications via the development of homotopy-theoretic approaches to higher category theory.  In particular, the framework of $(\infty, 1)$-categories, most prominently in the specific model of quasi-categories, have been used extensively in fields as diverse as topological quantum field theories and derived algebraic geometry.  But, reading the background literature can be daunting, especially if one tries to jump directly into the comprehensive treatments of Joyal \cite{joyal1}, \cite{joyal2} or Lurie \cite{lurie}.

There have been a number of attempts to make $(\infty,1)$-categories, and quasi-categories specifically, more accessible to non-experts, for example the expository papers of Antol\'in Camarena \cite{ac}, the author \cite{survey}, and Groth \cite{groth}; and in the recent book of Cisinski \cite{cisinski}.  But, those works often assume some prior familiarity with simplicial sets, or give the reader only a quick introduction to them.  The books of May \cite{may} and of Goerss and Jardine \cite{gj} both provide a comprehensive treatment of many important results about simplicial sets, although predating the current emphasis on quasi-categories, but again can be overwhelming for a novice.  

The purpose of this paper is to provide a still more gentle introduction to simplicial sets, taking a leisurely route to their definition.  As such, there is a good deal of overlap with the introductory paper of Friedman \cite{friedman}; some of the motivating ideas are also described more briefly by Dwyer \cite[\S 3]{dwyer}, but here we incorporate more of the category-theoretic perspective on the subject.  Our hope is that after reading this paper, one will be well-prepared to read and understand some of the other more introductory treatments of $(\infty,1)$-categories. 

We begin our exposition in Section \ref{simpcx} with a review of simplicial complexes, particularly how they can be viewed both geometrically and combinatorially.  The interplay between these two perspectives leads to the question of whether one can always model geometric or topological objects via purely combinatorial data.  As useful as they are for many purposes, we see in Sections \ref{orsimpcx} and \ref{ssets} that simplicial complexes do have some limitations when we try to apply certain constructions from topology.  Modifying their definition to solve these issues, we arrive at the definition of simplicial sets.

In Section \ref{cats} we introduce some basic ideas from category theory and use them to give a new definition of simplicial sets.  With that approach in hand, we revisit the question of how simplicial sets model topological spaces in Section \ref{ssetstop}, and then discuss the relationship with categories in Section \ref{ssetscat}.  Finally, in Section \ref{simpobj} we briefly describe some generalizations of simplicial sets and discuss why they are of interest.

\begin{thank}
	We would like to thank Walker Stern for assistance with some of the diagrams that appear here, and the referees for helpful comments.
\end{thank}

\section{Geometric and abstract simplicial complexes} \label{simpcx}

We begin our motivation for simplicial sets with a more broadly familiar concept: that of simplicial complexes.  In this section, we recall both the geometric and combinatorial approaches to defining simplicial complexes.  We begin with the geometric approach, and some definitions needed to define simplices.

Recall that points $v_0, v_1, \ldots, v_k$ in $\mathbb R^n$ are in \emph{general position} if any subset of them spans a strictly smaller hyperplane.  Given any set of points $v_0, \ldots , v_k \in \mathbb R^n$, their \emph{convex hull} is the smallest convex set containing them, namely, the set
\[ \left\{\sum_{i=0}^k \lambda_i v_i \mid \sum_{i=0}^k \lambda_i = 1, \lambda_i \in [0,1] \right\}. \]
Figure 1 illustrates the difference between the convex hull of points that are in general position compared to ones that are not.

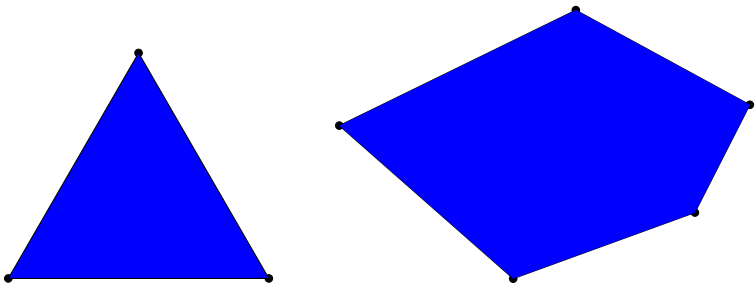
\begin{figure}[!htb]
	\begin{center}
		\begin{tabular}{c c c}
			\begin{tikzpicture}
				\foreach \tet/\lab in {-30/A,90/B,210/C}{
					\draw[fill=black] (\tet:2) circle (0.05);
					\path (\tet:2) node (\lab) {};
				}
				\path[fill=blue, opacity=0.3] (A.center) -- (B.center) -- (C.center) -- cycle;
				\draw (A.center) -- (B.center) -- (C.center) -- cycle;
			\end{tikzpicture} &  &
			\begin{tikzpicture}
				\foreach \tet/\len/\lab in {-30/2/A,10/2.5/B, 85/1.7/C,177/3/D,250/2/E,-90/0.5/F,99/1/G}{
					\draw[fill=black] (\tet:\len) circle (0.05); 
					\path (\tet:\len) node (\lab) {};
				};
				\draw (A.center) -- (B.center) -- (C.center) -- (D.center) -- (E.center) -- cycle;
				\path[fill=blue,opacity=0.3] (A.center) -- (B.center) -- (C.center) -- (D.center) -- (E.center) -- cycle;
			\end{tikzpicture}
		\end{tabular}
	\end{center}
	\caption{A simplex and a more general convex hull.}
\end{figure}

\begin{definition}
	Let $v_0, \ldots, v_k$ be $(k+1)$ points in general position in $\mathbb R^n$, for $n >k$.  The convex hull of these points is called a $k$-\emph{simplex}, denoted by $\Delta^n$.
\end{definition}

Examples of low-dimensional simplices are given in Figure 2.
\begin{figure}[!htb]
	\begin{center}
		\begin{tabular}{c c c c }
			$\Delta^0$ & $\Delta^1$ & $\Delta^2$ & $\Delta^3$ \\
			\begin{tikzpicture}
				\draw[fill=black] (0,0) circle (0.05);
			\end{tikzpicture} & \begin{tikzpicture}
				\draw[fill=black] (0,0) circle (0.05);
				\draw[fill=black] (2,0) circle (0.05);
				\draw (0,0) -- (2,0);
			\end{tikzpicture} &
			\begin{tikzpicture}[baseline=0.4cm]
				\foreach \tet/\lab in {-30/A,90/B,210/C}{
					\draw[fill=black] (\tet:1.5) circle (0.05);
					\path (\tet:1.5) node (\lab) {};
				}
				\path[fill=blue, opacity=0.3] (A.center) -- (B.center) -- (C.center) -- cycle;
				\draw (A.center) -- (B.center) -- (C.center) -- cycle;
			\end{tikzpicture} & 
			\begin{tikzpicture}[dot/.style={draw,circle,minimum size=1mm,inner sep=0pt,outer sep=0pt,fill=black}, scale=1.5, baseline=0.6cm]
				
				\coordinate[draw, dot] (3) at (0,{sqrt(2)},0);
				\coordinate[draw, dot] (2) at ({.5*sqrt(3)},0,-.5);
				\coordinate[draw, dot] (1) at (0,0,1);
				\coordinate[draw, dot] (0) at ({-.5*sqrt(3)},0,-.5);
				
				\path[fill=blue, opacity=0.3] (0) -- (1) -- (3) -- (2) -- cycle;
				
				\draw (0)--(1) (0)--(2) (0)--(3) (1)--(2) (1)--(3) (2)--(3);
				\path[fill=blue, opacity=0.3] (0.center) -- (1.center) -- (2.center) -- (3.center) -- cycle;
			\end{tikzpicture}
		\end{tabular}
	\end{center}
\caption{Examples of low-dimensional simplices.}
\end{figure}
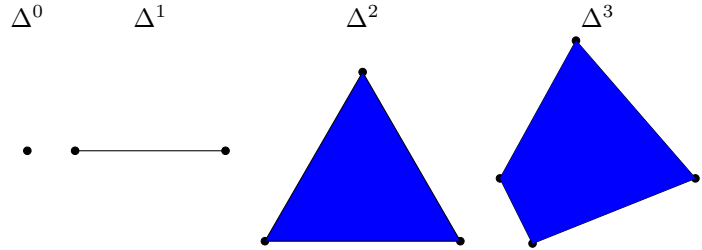
Since a 2-simplex is just a triangle, we often think of simplices (of any dimension) as ``generalized triangles."  The points $v_0, \ldots, v_k$ are called the \emph{vertices} of a simplex.  Restricting to a subset of the vertices, we get a simplex of a lower dimension, called a \emph{face} of the original simplex.  

Simplices are particularly nice convex spaces because they are completely determined by their vertices.  However, they are not topologically interesting, since they are always contractible.  We can obtain more interesting spaces by gluing simplices together in an appropriate way.

\begin{definition}
	A \emph{simplicial complex} $K$ in $\mathbb R^n$ is a finite set of simplices $P_1, \ldots ,P_m$ such that:
	\begin{enumerate}
		\item any face of $P_i$ is also in $K$, and
		
		\item for any $1 \leq i,j \leq m$, the intersection $P_i \cap P_j$ is either a face of both $P_i$ and $P_j$, or empty.
	\end{enumerate}
\end{definition}

An example is given in Figure 3.
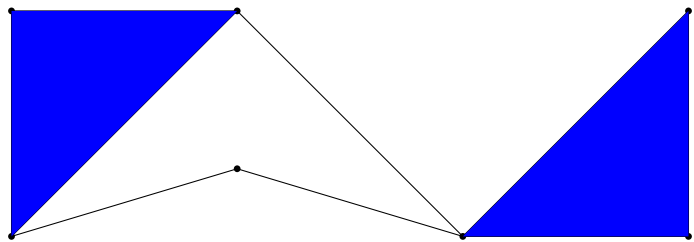
\begin{figure}[!htb]
	\begin{tikzpicture}[scale=3]
		\foreach \x/\y/\lab in {0/1/v1, 0/0/v2,1/1/v3,1/0.3/v4, 2/0/v5, 3/0/v6,3/1/v7}{
			\path (\x,\y) coordinate (\lab);
			\path[fill=black] (\x,\y) circle (0.015); 
			%		\path (\x,\y) node[label=above left:$\lab$] {};
		}
		\draw (v1)--(v2) (v1)--(v3) (v3)--(v2) (v4)--(v2) (v3)--(v5) (v4)--(v5) (v5)--(v6) (v5)--(v7) (v6)--(v7);
		\path[fill=blue, opacity=0.3] (v1)--(v2)--(v3)--(v1);
		\path[fill=blue, opacity=0.3] (v5)--(v6)--(v7)--(v5);
	\end{tikzpicture}
	\caption{A simplicial complex.}
\end{figure}

Many familiar objects can be thought of as simplicial complexes.  For example, polyhedra made up of triangles are simplicial complexes, such as the boundaries of the tetrahedron, the octahedron (as illustrated in Figure 4), and the icosahedron.  The solid tetrahedron is just a 3-simplex.  

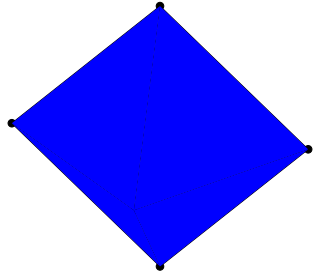
\begin{figure}[!htb]
\tdplotsetmaincoords{60}{100}
\begin{tikzpicture}[dot/.style={draw,circle,minimum size=1mm,inner sep=0pt,outer sep=0pt,fill=black}, scale=2, baseline=0.6cm,tdplot_main_coords]
	
	\coordinate[draw, dot] (T) at (0,0,1);
	\coordinate[draw, dot] (B) at (0,0,-1);
	\coordinate[draw, dot] (L) at (-1,0,0);
	\coordinate[draw, dot] (R) at (1,0,0);
	\coordinate[draw, dot] (F) at (0,1,0);
	\coordinate[draw, dot] (H) at (0,-1,0);
	\draw (R.center)-- (H.center)--(L.center) --(F.center) ;
	\draw (T.center) -- (H.center);
	\draw (T.center) -- (L.center);
	\draw (T.center) -- (F.center);
	
	\draw[shorten <= 4mm,shorten >=4mm, line width=2mm,white] (T.south) -- (R.north);
	\draw(T.center) -- (R.center);
	
	\draw (B.center) -- (R.center);
	\draw (B.center) -- (F.center);
	\draw (B.center) -- (H.center); 
	\draw (B.center) -- (L.center); 
	
	\draw[shorten <= 4mm,shorten >=4mm,line width=2mm,white] (F) -- (R);
	\draw (F.center) -- (R.center);
	
	\path[fill=blue,opacity=0.3] (T.center)-- (R.center) --(F.center) -- cycle;
	\path[fill=blue,opacity=0.2] (T.center) --(R.center) -- (H.center)--cycle;
	\path[fill=blue,opacity=0.25] (B.center) --(R.center) -- (H.center)--cycle;
	\path[fill=blue,opacity=0.4] (B.center)-- (R.center) --(F.center) -- cycle;
\end{tikzpicture}
\caption{The octahedron as a simplicial complex.}
\end{figure}

As we have already observed, simplices are completely determined by their vertices, and with a bit more care, we can also describe simplicial complexes purely combinatorially, as follows.

\begin{definition}
	An \emph{abstract simplicial complex} $K$ is a pair $(V_K, S_K)$, where $V_K$ is a finite set (whose elements are called \emph{vertices}) and $S_K$ is a set of nonempty finite subsets of $V_K$ (called \emph{simplices}) such that all singleton subsets of $V_K$ are in $S_K$, and, if $\sigma \in S_K$ and $\sigma' \subseteq \sigma$, then $\sigma' \in S_K$.
\end{definition}

Here, we are including the adjective ``abstract" to emphasize the difference from the geometric definition, but we claim that the two structures are equivalent.  We have essentially described how to extract an abstract simplicial complex from a geometric one.

\begin{example} \label{scxex}
	An example of an abstract simplicial complex $K$ is given by $V_K = \{a,b,c,d\}$ and
	\[ S_K=\{\{a\}, \{b\}, \{c\}, \{d\}, \{a,b\}, \{a,c\}, \{b,c\}, \{b,d\}, \{c,d\}, \{a,c,d\} \}. \]
	We leave verification of the conditions to the reader.
\end{example}

\begin{example}
	On the other hand, if we define $L$ by $V_L = \{a,b,c\}$ and
	\[ S_L = \{\{a\}, \{b\}, \{c\}, \{a,b\}, \{a,c\}, \{a,b,c\}\}, \]
	we do not obtain an abstract simplicial complex, since the 2-simplex $\{a,b,c\}$ must have face $\{b,c\}$ included in $S_L$.  
\end{example}

In the other direction we can \emph{geometrically realize} an abstract simplicial complex to a geometric one by assigning actual points in $\mathbb R^n$, for sufficiently large $n$, to the elements of $V_K$ and including $k$-simplices as specified by the $(k+1)$-element sets in $S_K$.  For small examples, we can just draw them; for instance, the simplicial complex in Example \ref{scxex} can be realized as:
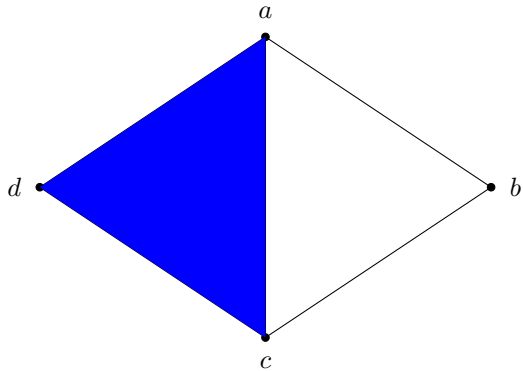
\[	\begin{tikzpicture}
			\foreach \x/\y/\lab/\pos in {0/2/a/above,0/-2/c/below,3/0/b/right,-3/0/d/left}{
				\draw[fill=black] (\x,\y) circle (0.05);
				\path (\x,\y) node[label=\pos:$\lab$] (\lab) {};
			};
			\draw (a.center) -- (b.center) -- (c.center) -- (d.center) -- cycle;
			\draw (a.center) -- (c.center); 
			\path[fill=blue,opacity=0.3] (a.center) -- (c.center) -- (d.center) -- cycle;
		\end{tikzpicture}. \]

However, we need a more formal definition, especially for working with larger examples.  To start, we want to be sure that $n$ is sufficiently large, in particular at least as large as the cardinality of the largest set in $S_K$.  One approach is the following construction.  Let $n=|V_K|$, and let $e_i$ for $1 \leq i \leq n$ be the standard basis vectors in $\mathbb R^n$.  Specify a bijection $\varphi \colon V_K \rightarrow \{1, \ldots n\}$ so that each vertex has a distinct label.  Then the geometric realization of $K$ is given by
\[ |K| = \bigcup_{\sigma \in S_K} c(\sigma) \] 
where $c(\sigma)$ denotes the convex hull of the set
\[ \{e_{\varphi(s)} \}_{s\in \sigma}. \]
See Figure 5 for a depiction of the output of this process for the 2-simplex.
\begin{figure}
	\tdplotsetmaincoords{100}{120}
	\begin{tikzpicture}[tdplot_main_coords]
		\draw[->] (-1,0,0) -- (3,0,0) node[above] {$x$};
		\draw[->] (0,-1,0) -- (0,3,0) node[above] {$y$};
		\draw[->] (0,0,-1) -- (0,0,3) node[above] {$z$};
		\path (2,0,0) coordinate (y);
		\path (0,0,2) coordinate (x);
		\path (0,2,0) coordinate (z);
		\draw[fill=black] (y) circle (0.05); 
		\draw[fill=black] (x) circle (0.05); 
		\draw[fill=black] (z) circle (0.05);
		\draw[thick] (x) -- (y) -- (z)-- (x);
		\path[fill=blue,opacity=0.3]  (x) -- (y) -- (z)-- (x);
	\end{tikzpicture}
	\caption{Standard geometric realization of $\Delta^2$.}
\end{figure}
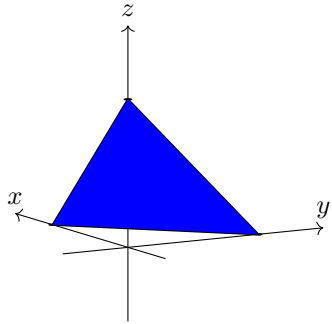

Observe that there is an overlap in notation here, with $|V_K|$ denoting the cardinality of the set $|V_K|$, and $|K|$ denoting the geometric realization of the simplicial complex $K$.  Both usages are standard, so at the risk of confusion we continue to use them here.

If we take the geometric realization of an abstract simplicial complex, and then describe that geometric object combinatorially, we recover the original abstract simplicial complex.  Working the other way around, if we start with a geometric simplicial complex, describe it combinatorially, and then apply geometric realization, we might not get the identical simplicial complex but we do get one that homeomorphic to it.  

\section{From simplicial complexes to oriented simplicial complexes} \label{orsimpcx}

We have now seen a close relationship between simplicial complexes and at least some kinds of topological spaces. Given an abstract simplicial complex, we can take its geometric realization to get a topological space.  If we move away from working linearly, we can triangulate topological spaces to get simplicial complexes.  Doing so helps us to compute invariants such as the Euler characteristic, homology groups, or the fundamental group for a given topological space.  

Now that we have discussed the similarities between abstract simplicial complexes and topological spaces that can be realized as simplicial complexes, we turn our attention to some differences, namely ways in which the two kinds of data are not compatible with one another. 

Let us begin with the product operation on simplicial complexes.  Here, we take the perspective that a product should be given by a particular universal property.  Specifically, if $K$ and $L$ are simplicial complexes, then their product $K \times L$ should have the property that, given any other simplicial complex $M$ equipped with maps $M \rightarrow K$ and $M \rightarrow L$, there is a unique map $M \rightarrow K \times L$ making the diagram
\[ \xymatrix{M \ar[dr] \ar[ddr] \ar[drr] && \\
	& K \times L \ar[r] \ar[d] & L \\
	& L & } \]
commute. 

We claim that the following construction satisfies this property; we leave the verification to the reader.  Let $V_{K \times L} = V_K \times V_L$, so the vertices of the product are given pairs of vertices of the original simplicial complexes.  To define the simplices of $K \times L$, consider the projection maps $p_1 \colon V_K \times V_L \rightarrow V_K$, given by $(v, v') \mapsto v$, and $p_2 \colon V_K \times V_L \rightarrow V_L$, given by $(v,v') \mapsto v'$. Now, we define
\[ S_{K \times L} =\{\text{all subsets } \sigma \text{ of } V_{K \times L} \text{ with } p_1(\sigma) \in S_K, p_2(\sigma) \in S_L \}. \]

\begin{example}
	Let $\Delta^1$ be the 1-simplex and $\partial \Delta^1$ be its boundary, which consists of two points.  Thus, $V_{\Delta^1} = \{v_0, v_1\}$ and $V_{\partial \Delta^1}= \{w_0, w_1\}$, with
	\[ S_{\Delta^1}=\{\{v_0\}, \{v_1\}, \{v_0, v_1\}\} \] and
	\[ S_{\partial \Delta^1}= \{\{w_0\}, \{w_1\}\}. \]
	Using the definition of product described above, we get
	\[ V_{\Delta^1 \times \partial \Delta^1} = \{(v_0, w_0), (v_0, w_1), (v_1, w_0), (v_1, w_1)\} \] and
	\[ S_{\Delta^1 \times \partial \Delta^1} = \{\{(v_0, w_0)\}, \{(v_0, w_1)\}, \{(v_1, w_0)\}, \{(v_1, w_1)\}, \{(v_0, w_0), (v_1, w_0)\}, \{(v_0, w_1), (v_1, w_1)\} \}. \]
	Thus, the geometric realization $|\Delta^1 \times \partial \Delta^1|$ looks as we would expect it to: 
	
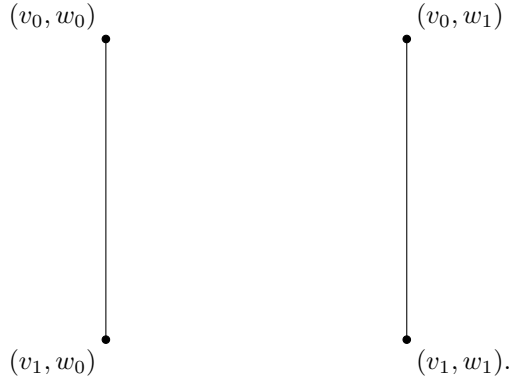
\begin{center}
		\begin{tikzpicture}[dot/.style={draw,circle,minimum size=1mm,inner sep=0pt,outer sep=0pt,fill=black},scale=2]
			\coordinate[dot] (A) at (0,1);
			\coordinate[dot] (B) at (0,-1);
			\coordinate[dot] (C) at (2,1); 
			\coordinate[dot] (D) at (2,-1);
			\node at (A) [above left] {$(v_0,w_0)$};
			\node at (B) [below left] {$(v_1,w_0)$};
			\node at (C) [above right] {$(v_0,w_1)$};
			\node at (D) [below right] {$(v_1,w_1)$.};
			\draw (A.center) -- (B.center);
			\draw (C.center) -- (D.center);
		\end{tikzpicture}
	\end{center}
\end{example}

However, it is not always the case that the geometric realization of a simplicial complex is the same as the product of the geometric realizations.

\begin{example}
	Consider $\Delta^1 \times \Delta^1$.  Using notation from above, we get
	\[ V_{\Delta^1 \times \Delta^1}= \{(v_0, v_0), (v_0, v_1), (v_1, v_0), (v_1, v_1)\}. \]  
	One can check that, using the criterion above, the elements of $S_{\Delta^1 \times \Delta^1}$ are all the subsets of $V_{\Delta^1 \times \Delta^1}$.  In particular, $\Delta^1 \times \Delta^1$ is a 3-simplex:
	
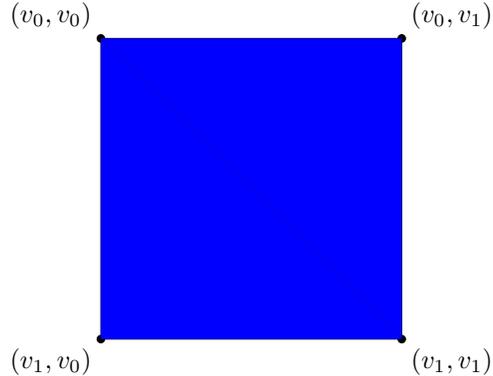
\begin{center}
		\begin{tikzpicture}[dot/.style={draw,circle,minimum size=1mm,inner sep=0pt,outer sep=0pt,fill=black},scale=2]
			\coordinate[dot] (A) at (0,1);
			\coordinate[dot] (B) at (0,-1);
			\coordinate[dot] (C) at (2,1); 
			\coordinate[dot] (D) at (2,-1);
			\node at (A) [above left] {$(v_0,v_0)$};
			\node at (B) [below left] {$(v_1,v_0)$};
			\node at (C) [above right] {$(v_0,v_1)$};
			\node at (D) [below right] {$(v_1,v_1)$.};
			\draw (A.center) -- (B.center);
			\draw (C.center) -- (D.center);
			
			\draw (B.center) -- (C.center);
			\draw[shorten <= 4mm,shorten >=4mm, line width=3mm,white] (A.center) -- (D.center);
			\draw (A.center) -- (D.center);
			\draw (A.center) -- (C.center);
			\draw (B.center) -- (D.center);
			\path[fill=blue,opacity=0.2] (A.center) --(B.center) -- (D.center) -- cycle;
			\path[fill=blue,opacity=0.3] (A.center) --(D.center) -- (C.center) -- cycle;
		\end{tikzpicture}
	\end{center}
	However, geometrically we expect $\Delta^1 \times \Delta^1$ to be a square, and in particular only 2-dimensional.
\end{example}

So, the geometric realization and product operations do not commute with one another, and in particular simplicial complexes do not model topological spaces very well if we want to take products.  We thus ask whether we can add some additional structure to simplicial complexes to remedy this difficulty.  The answer lies in the following definition.

\begin{definition}
	An \emph{oriented} simplicial complex is a simplicial complex with a partial ordering on its vertices.
\end{definition}

Given any abstract simplicial complex, we can make it an oriented simplicial complex by totally ordering its vertices.  Then each 1-simplex has a ``direction" given by the ordering.  Geometrically, instead of thinking of a 1-simplex as a segment
\[ \begin{tikzpicture}
		\draw[fill=black] (0,0) circle (0.05);
		\draw[fill=black] (2,0) circle (0.05);
		\draw (0,0) -- (2,0);
\end{tikzpicture} \]
we instead think of it as an arrow
\[ \begin{tikzpicture}
		\draw[fill=black] (0,0) circle (0.05);
		\draw[fill=black] (2,0) circle (0.05);
		\path (0,0) node (A) {};
		\path (2,0) node (B) {};
		\draw[->] (A) to (B);
	\end{tikzpicture} \]
indicating the ordering.  Thus the oriented 2-simplex 
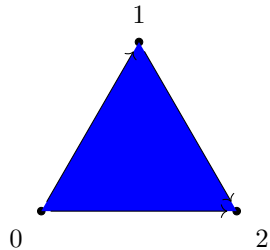
\[	\begin{tikzpicture}
	\foreach \tet/\lab/\nom/\pos in {-30/A/2/below right,90/B/1/above,210/C/0/below left}{
		\draw[fill=black] (\tet:1.5) circle (0.05);
		\path (\tet:1.5) node[label=\pos:$\nom$] (\lab) {};
	}
	\path[fill=blue, opacity=0.3] (A.center) -- (B.center) -- (C.center) -- cycle;
	\draw[->] (C) to (B);
	\draw[->] (B) to (A);
	\draw[->] (C) to (A);
\end{tikzpicture} \]
has a particular arrangement of arrows following from the ordering on the vertices.  Observe that we have labeled the vertices here by $0,1,2$ to emphasize their ordering.  We follow this convention for the oriented $n$-simplex for any $n \geq 0$ and denote it by $\Delta[n]$ to distinguish it from its non-oriented counterpart.  We can geometrically realize an oriented simplicial complex just as before; we simply forget the orientation.  Thus, the geometric realization of $\Delta[n]$ is $\Delta^n$, where we implicitly take the latter to be a geometric object rather than its combinatorial counterpart.

\begin{example} \label{ssetexs}
	Let us introduce some other important examples.  We have already introduced the $n$-simplex $\Delta[n]$ for any $n \geq 0$.  We can also consider its \emph{boundary} $\partial \Delta[n]$, which is the oriented simplicial complex obtained by removing the single top-dimensional simplex from it.  We thus obtain the union of all the lower-dimensional faces of $\Delta[n]$.  So, for example, $\partial \Delta[1]$ consists of two points, just as in the unoriented case above, and the $\partial \Delta[2]$ is given by 
	
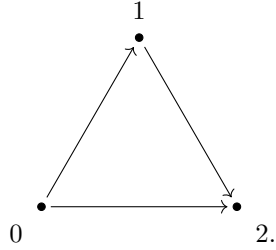
\[ \begin{tikzpicture}
		\foreach \tet/\lab/\nom/\pos in {-30/A/2./below right,90/B/1/above,210/C/0/below left}{
			\draw[fill=black] (\tet:1.5) circle (0.05);
			\path (\tet:1.5) node[label=\pos:$\nom$] (\lab) {};
		}
		
		\draw[->] (C) to (B);
		\draw[->] (B) to (A);
		\draw[->] (C) to (A);
	\end{tikzpicture} \]  
	Observe that $\partial \Delta[0] = \varnothing$. 
	
	We can also remove one of the highest-dimensional faces of $\partial \Delta[n]$ when $n \geq 1$ and the ordering on the vertices enables us to specify precisely which face is removed.  Using the labeling of the vertices of $\partial \Delta[n]$ by $0, 1, \ldots, n$, then for any $0 \leq k \leq n$, we denote by $\Lambda^k[n]$ the oriented simplicial complex obtained from $\partial \Delta[n]$ by removing the unique $(n-1)$-dimensonal face not containing the vertex labeled by $k$.  This simplicial complex is called a \emph{horn}, as suggested by the pictures in Figure 6.
	\begin{figure}[!htb]
	\begin{tabular}{ c c c }
		\begin{tikzpicture}[scale=0.7]
			\foreach \tet/\lab/\nom/\pos in {-30/A/2/below right,90/B/1/above,210/C/0/below left}{
				\draw[fill=black] (\tet:1.5) circle (0.07);
				\path (\tet:1.5) node[label=\pos:$\nom$] (\lab) {};
			}
			\draw[->] (C) to (B);
			\draw[->] (C) to (A);
		\end{tikzpicture}
		& $\quad$
		\begin{tikzpicture}[scale=0.7]
			\foreach \tet/\lab/\nom/\pos in {-30/A/2/below right,90/B/1/above,210/C/0/below left}{
				\draw[fill=black] (\tet:1.5) circle (0.07);
				\path (\tet:1.5) node[label=\pos:$\nom$] (\lab) {};
			}
			\draw[->] (C) to (B);
			\draw[->] (B) to (A);
		\end{tikzpicture}
		& $\quad$
		\begin{tikzpicture}[scale=0.7]
			\foreach \tet/\lab/\nom/\pos in {-30/A/2/below right,90/B/1/above,210/C/0/below left}{
				\draw[fill=black] (\tet:1.5) circle (0.07);
				\path (\tet:1.5) node[label=\pos:$\nom$] (\lab) {};
			}
			\draw[->] (B) to (A);
			\draw[->] (C) to (A);
		\end{tikzpicture}
	\end{tabular}
		\caption{The horns $\Lambda^0[2]$, $\Lambda^1[2]$, and $\Lambda^2[2]$.}
		\end{figure}
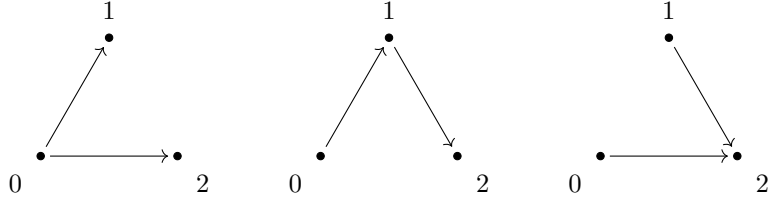
\end{example}

Given oriented simplicial complexes $K$ and $L$, we define their product $K \times L$ to have vertices $V_{K \times L} = V_K \times V_L$.  Observe that this set is still partially ordered.  Its simplices are given by
\[ S_{K \times L}=\{\text{totally ordered subsets } \sigma \text{ of } V_{K \times L} \text{ with } p_1(\sigma) \in S_K, p_2(\sigma) \in S_L\}. \]
Again, we leave it to the reader to verify that this definition has the desired universal property.

Let us now revisit the problematic example from before.

\begin{example}
	Let us take the product of oriented simplicial complexes $\Delta[1] \times \Delta[1]$.  We label the vertices of $\Delta[1]$ by $v_0$ and $v_1$ as before, and order them by $v_0 < v_1$.  Then the set of vertices in the product is the same as before.  However, now we do not get all possible subsets as simplices, because two of the vertices, $(v_0, v_1)$ and $(v_1, v_0)$, are incomparable.  The resulting oriented simplicial complex does have geometric realization a square, just as we expect topologically.
\end{example}

\section{From oriented simplicial complexes to simplicial sets} \label{ssets}

By incorporating the structure of an orientation, we have modified simplicial complexes so that they better model topological spaces.  However, they still do not have all the properties that we want.  An additional problem arises when we take quotients.

For example, take an interval and glue its ends together to get a circle:
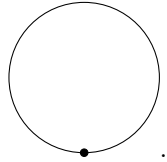
\[ \begin{tikzpicture}
	\draw (0,0) circle (1);
	\draw[fill=black] (-90:1) circle (0.05);
\end{tikzpicture}. \]  
Using the definition of quotient spaces, this procedure is well-defined.  What if we try to do this construction on the level of (oriented) simplicial complexes?  We could start with $\Delta[1]$, and then identify its vertices.  But, we have assumed that simplices are completely determined by their vertices, so a simplicial complex with one vertex can only be the 0-simplex $\Delta[0]$.  

We can still find an (oriented) simplicial complex whose geometric realization is homeomorphic to a circle, for example, the boundary of the 2-simplex.  But, this simplicial complex has three vertices and three 1-simplices, and there does not seem to be a natural way to get such an object as a quotient of a 1-simplex, especially since there are many other possibilities, for example the boundary of any regular polygon, such as
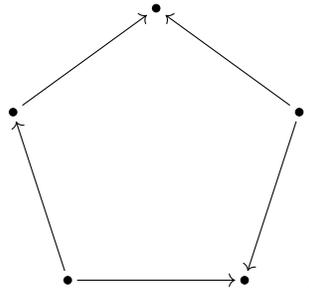
\[ \begin{tikzpicture}
	\foreach \tet/\lab in {-54/A,18/B,90/C,162/D,234/E}{
		\draw[fill=black] (\tet:2) circle (0.05);
		\path (\tet:2) node (\lab) {};
	};
	\draw[->] (E) -- (A); 
	\draw[->] (E) -- (D);
	\draw[->] (D) -- (C);
	\draw[->] (B) -- (C);
	\draw[->] (B) to (A);
\end{tikzpicture}. \]
But, intuitively what we really want is something like
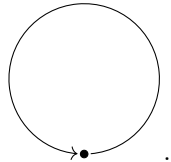
\[ \begin{tikzpicture}
	\draw[->] (-85:1) arc (-85:265:1);
	\draw[fill=black] (-90:1) circle (0.05);
\end{tikzpicture}. \]

We claim that the solution here is to abandon the condition that the simplices in a simplicial complex have to be determined by their vertices. In other words, we want to consider more general ``simplices" that no longer look triangular.  

To make this transition, let us think about simplicial complexes a bit differently.  Previously, given a simplicial complex $K$, we put all simplices of all dimensions together into a single set $S_K$, and the dimension of each simplex was implicit from the number of vertices specifying it.  If simplices are no longer to be determined by their vertices, this method no longer works, and so we need extra data to specify the dimension of a given simplex.  We encode this data via \emph{face maps}, which take an $n$-simplex to one of its $(n-1)$-simplex faces.  For example, given a 1-simplex
\[ \bullet_{v_0} \overset{e}{\longrightarrow} \bullet_{v_1}, \]
we have two face maps picking out the two vertices $v_0$ and $v_1$.   In the more general context, we might have another 1-simplex $f$ with these same vertices, which was not allowed previously.

Let us describe this structure in detail for oriented simplicial complexes, and then use it to define the more general structures we want.  For an oriented simplicial complex $K$, let $K_n$ denote its set of $n$-simplices.  Note that $K_0$ coincides with the set $V_K$ of vertices of $K$.

There are $n+1$ face maps $K_n \rightarrow K_{n-1}$, given by taking the face opposite each of the vertices; we can use the ordering on the vertices to specify which is which.  More specifically, let us consider an $n$-simplex with vertices labeled by $v_0, v_1, \ldots, v_n$.  Given any $0 \leq i \leq n$, we define the $i$-\emph{face} to be the $(n-1)$-simplex with vertices $v_0, \ldots, v_{i-1}, v_{i+1}, \ldots, v_n$.  Then for any $0 \leq i \leq n$, the $i$-th \emph{face map} $d_i \colon K_n \rightarrow K_{n-1}$ takes any $n$-simplex to its $i$-face.

For higher-dimensional simplices, observe that there are some relations when we take iterated face maps. For example, consider $\Delta[2]$:
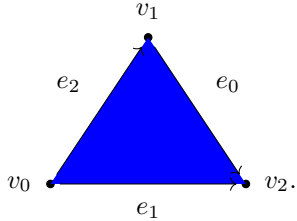
\[ \begin{tikzpicture}
	\path (210:1.5) node[label=left: $v_0$] (C) {}; 
	\path (90:1.2) node[label=above:$v_1$] (B) {};
	\path (-30:1.5) node[label=right:$v_2$.] (A) {};
	\draw[fill=black] (A.center) circle (0.05);
	\draw[fill=black] (B.center) circle (0.05);
	\draw[fill=black] (C.center) circle (0.05);
	\path (0,0) node {$\sigma$};
	\path[fill=blue,opacity=0.3] (A.center) -- (B.center) -- (C.center) -- cycle;
	\draw[->] (C) to node[label=above left:$e_2$]{} (B);
	\draw[->] (B) to node[label=above right:$e_0$]{} (A);
	
	\draw[->] (C) to node[label=below:$e_1$]{} (A);
\end{tikzpicture} \]
Applying each face map to the 2-simplex $\sigma$, we get 
\[ d_0(\sigma) = e_0, d_1(\sigma)=e_1, d_2(\sigma) = e_2. \]
(Note that we have chosen the names of the 1-simplices strategically here!)  Taking faces of these 1-simplices, we get
\[ \begin{aligned}
	d_0d_0(\sigma) = d_0(e_0)=v_2, && d_1d_0(\sigma) = d_1(\sigma)=v_1, \\
	d_0d_1(\sigma) = d_0(e_1)=v_2, && d_1d_1(\sigma) = d_1(e_1) = v_0, \\
	d_0d_2(\sigma)=d_0(e_2) = v_1, && d_1d_2(\sigma) = d_1(e_2)=v_0.
\end{aligned} \]
We thus observe the relations
\[ \begin{aligned}
	d_0d_0 & = d_0d_1 \\
	d_0d_2 & = d_1d_0 \\
	d_1d_2 & = d_1d_1. 
	\end{aligned} \]
One can check that similar relations hold for higher face maps, satisfying the general condition that $d_id_j = d_{j-1}d_i$ for $i<j$.
	
Formalizing this idea, we make the following definition.

\begin{definition}
	A $\Delta$-\emph{set} $K$ consists of sets $K_n$ for each $n \geq 0$, together with, for each $n \geq 1$, \emph{face maps} $d_i \colon K_n \rightarrow K_{n-1}$, satisfying $d_id_j=d_{j-1}d_i$ for $i<j$.
	
	A map of $\Delta$-sets $f \colon K \rightarrow L$ consists of functions $f_n \colon K_n \rightarrow L_n$ for all $n \geq 0$ that commute with face maps, so that $d_i f_n = f_{n-1} d_i$ for all $n \geq 1$ and $0 \leq i \leq n$.
\end{definition}

We can now describe the circle with a single vertex as a $\Delta$-set $K$ with $K_0$ and $K_1$ each consisting of a single point, and $K_n = \varnothing$ for $n \geq 2$.  The two face maps $K_1 \rightarrow K_0$ are the only thing they can be.  In this case, we also have a quotient map $\Delta[1] \rightarrow K$ given by collapsing both vertices of $\Delta[1]$ to the single vertex of $K$.

Topologically, such quotient maps always exist, but between $\Delta$-sets they need not, particularly when we collapse higher-dimensional simplices down to lower-dimensional ones.  For example, consider the quotient of the 2-simplex $\Delta[2]$ by identifying the edge opposite the vertex $0$.  The resulting $\Delta$-set $D$ has two 0-simplices, two 1-simplices, and a single 2-simplex.  However, there is no map of $\Delta$-sets $\Delta[2] \rightarrow D$, since there is no 1-simplex in $D$ to which we can send the 1-simplex that was collapsed.

The new requirement that maps preserve dimension presents an obstacle to defining such quotient maps.  What we really want is a way to think of any $n$-simplex as an ``honorary" $k$-simplex for every $k>n$.  In our example above, we can thus send the 1-simplex to the 0-simplex it gets collapsed to by thinking of it as a 1-simplex in a suitable way.

How can we formalize this idea?  We want to include ``upward maps" $s_i \colon K_n \rightarrow K_{n+1}$ whose images consist of these degenerate $(n+1)$-simplices.  But this process becomes complicated as we try to fit such maps in with face maps.  Let us look at an example.

When $n=0$, we think of a 0-simplex as a degenerate simplex in a straightforward way: we think of it as a ``very small" 1-simplex between two copies of the vertex itself.  The two face maps from this degenerate 1-simplex are completely determined.

However, even when $n=1$ we find we need to make choices.  Consider a 1-simplex
\[ v_0 \overset{e}{\longrightarrow} v_1. \]
How can we think of it as a degenerate 2-simplex?  We can consider this 1-simplex as comprising two of the faces of the 2-simplex, and then having the third face be a degenerate 1-simplex.  But there are two choices: we can take a degenerate edge on $v_0$, or on $v_1$.  Thus, we have two natural maps $K_1 \rightarrow K_2$, and we denote them by $s_0$ and $s_1$, the subscript denoting which vertex is being regarded as a degenerate edge.  See Figure 7 for a depiction of these maps.
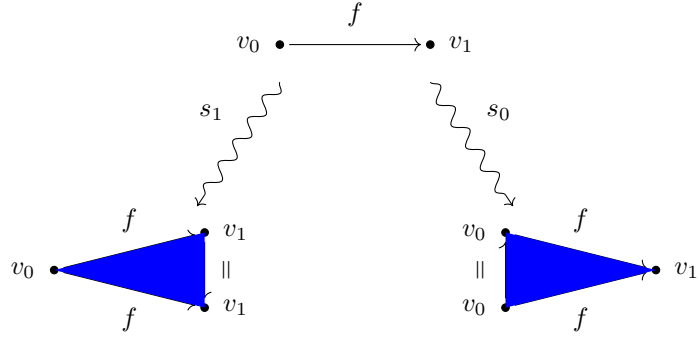
\begin{figure}[!htp]
	\begin{tikzpicture}
		\path (0,0) node[label=left:$v_0$] (v0) {};
		\path (2,0) node[label=right:$v_1$] (v1) {};
		\draw[fill=black] (v0.center) circle (0.05);
		\draw[fill=black] (v1.center) circle (0.05);
		\draw[->] (v0) to node[label=above:$f$]{} (v1);
		
		\path (-3,-3) node[label=left:$v_0$] (vl0) {};
		\path (-1,-3.5) node[label=right:$v_1$] (vlb1) {};
		\path (-1,-2.5) node[label=right:$v_1$] (vlt1) {};
		\draw[fill=black] (vl0.center) circle (0.05);
		\draw[fill=black] (vlb1.center) circle (0.05);
		\draw[fill=black] (vlt1.center) circle (0.05);
		\draw[->] (vl0) to node[label=above:$f$]{} (vlt1);
		\draw[->] (vl0) to node[label=below:$f$]{} (vlb1);
		\draw[->] (vlt1) to node[label=right:$\rotatebox{90}{=}$]{}  (vlb1);
		\path (-1.6,-3) node {$\rotatebox{90}{=}$};
		\path[fill=blue, opacity=0.2] (vl0.center) -- (vlt1.center) -- (vlb1.center) -- cycle;
		
		\path (5,-3) node[label=right:$v_1$] (vr0) {};
		\path (3,-3.5) node[label=left:$v_0$] (vrb1) {};
		\path (3,-2.5) node[label=left:$v_0$] (vrt1) {};
		\draw[fill=black] (vrb1.center) circle (0.05) (vrt1.center) circle (0.05) (vr0.center) circle (0.05);
		\draw[->] (vrb1) to node[label=left:$\rotatebox{90}{=}$]{} (vrt1);
		\draw[->] (vrt1) to node[label=above:$f$]{} (vr0);
		\draw[->] (vrb1) to node[label=below:$f$] {} (vr0);
		\path (3.6,-3) node {$\rotatebox{90}{=}$};
		\path[fill=blue,opacity=0.2] (vrb1.center) -- (vrt1.center) -- (vr0.center) -- cycle;
		
		\begin{scope}[decoration=snake]
			\draw[decorate] (0,-0.5) to node[label=above left:$s_1$]{} (-1,-2);
			\draw[->] (-1,-2) to (-1.1,-2.15);
		\end{scope}
		
		\begin{scope}[decoration=snake]
			\draw[decorate] (2,-0.5) to node[label=above right:$s_0$]{} (3,-2);
			\draw[->] (3,-2) to (3.1,-2.15);
		\end{scope}
	\end{tikzpicture}
\caption{The degeneracy maps $s_0$ and $s_1$}
\end{figure}

A similar argument shows that there should be $n+1$ degeneracy maps $K_n \rightarrow K_{n+1}$ for every $n \geq 0$.  We can think of each as ``doubling" one of the vertices and including a degenerate 1-simplex.  Before proceeding to the definition below, we invite the reader to investigate some relations between the face and degeneracy maps.

We can now state the definition of a simplicial set.

\begin{definition} \label{simpsetdef}
	A \emph{simplicial set} $X$ is a collection of sets $X_0, X_1, X_2, \ldots$ together with face maps $d_i \colon X_n \rightarrow X_{n-1}$ for $0 \leq i \leq n$ and degeneracy maps $s_i \colon X_n \rightarrow X_{n+1}$ for $0 \leq i \leq n$, satisfying the relations
	\[ \begin{aligned}
	d_id_j & = d_{j-1}d_i & i<j \\
	d_is_j & = s_{j-1}d_i & i<j \\
	d_js_j & =d_{j+1}s_j =\id & \\
	d_i s_j & = s_jd_{i-1} & i>j+1 \\
	s_is_j & s_{j+1}s_i & i \leq j.
	\end{aligned} \]
\end{definition}

Thus, a simplicial set looks like a diagram of sets and functions
\[ \xymatrix@1{X_0 \ar@<-.5ex>[r] & X_1 \ar[l] \ar@<-.5ex>[l] \ar@<-.5ex>[r] \ar@<-1ex>[r] & X_2 \ar[l] \ar@<-.5ex>[l] \ar@<-1ex>[l] \ar@<-.5ex>[r] \ar@<-1ex>[r] \ar@<-1.5ex>[r] & \cdots. \ar[l] \ar@<-.5ex>[l] \ar@<-1ex>[l] \ar@<-1.5ex>[l]} \]

One significant difference between $\Delta$-sets and simplicial sets is that the latter, unless empty, inherently have infinitely many simplices.  For example, as a $\Delta$-set, $\Delta[1]$ has a single 1-simplex and no simplices of higher dimensions. As a simplicial set, it has degenerate simplices of arbitrarily high dimensions.  Other examples of simplicial complexes and $\Delta$-sets that we have considered can be reinterpreted as simplicial sets, for example the boundaries $\partial \Delta[n]$ and horns $\Lambda^k[n]$.

We claim that simplicial sets are the combinatorial model for topological spaces that we have been looking for.  To formalize what mean by this statement, namely, to establish some kind of equivalence between simplicial sets and topological spaces, we need to introduce some of the language of category theory.  Along the way, we can also give a more efficient definition of simplicial sets.

\section{Categories, functors, and a new approach} \label{cats}

To talk about simplicial sets in a more systematic way, it helps to describe them equivalently in a way that uses the language of categories and functors.  We begin by recalling the definition of a category.

\begin{definition}
	A \emph{category} $\mathcal C$ consists of a collection of \emph{objects} together with:
	\begin{enumerate}
		\item for any ordered pair $(a,b)$ of objects, a set $\Hom_\mathcal C(a,b)$ of \emph{morphisms} $a \rightarrow b$, and
		
		\item for any ordered triple $(a,b,c)$ of objects, a function \[ \Hom_\mathcal C(a,b) \times \Hom_\mathcal C(b,c) \rightarrow \Hom_\mathcal C(a,c), \] 
		giving composites of morphisms, say $(f,g) \mapsto g \circ f$, satisfying:
		
		\begin{itemize}
			\item composition is associative, so, if we have
			\[ \xymatrix@1{a \ar[r]^f & b \ar[r]^g & c \ar[r]^h & d}, \] then $h \circ(g \circ f) = (h \circ g) \circ f$; and
			
			\item for each object $b$ there is an identity morphism $\id_b \in \Hom_\mathcal C(b,b)$ such that, if $f \colon a \rightarrow b$, then $\id_b \circ f = f$, and, if $g \colon b \rightarrow c$, $g \circ \id_b = g$.
		\end{itemize}
	\end{enumerate}
\end{definition}

Observe that we have been somewhat ambiguous in describing the ``collection" of objects of a category, whereas we have specified that each $\Hom_\mathcal C(x,y)$ is a set.  Categories whose objects form a set are called \emph{small}; if the objects form a proper class, we sometimes call such a category \emph{large} for emphasis.   

Many familar objects in mathematics, together with appropriate functions between them, form large categories.  Let us look at a few examples.

\begin{example}
	The class of all sets and all functions between them forms a category we denote by $\Sets$.  
\end{example}

\begin{example}
	The category $\Top$ has objects topological spaces and morphisms continuous maps between them.
\end{example}

\begin{example}
	The category $\mathcal{ASC}$ has objects abstract simplicial complexes.  Functions between simplicial complexes are given by set functions on vertices, together with compatible functions between simplices.  Here, compatibility means not only agreeing with the functions on vertices, but also respecting face maps.
\end{example}

We would like to describe a category of simplicial sets similarly, so that a morphism is given by functions between $n$-simplices for each $n \geq 0$, commuting with all face and degeneracy maps.  To make this description more precise, we turn to an alternative description of simplicial sets.  To do so, we first need the following example of a small category.

\begin{example}
	For each $n \geq 0$, consider the finite ordered set
	\[ [n] = \{0 \leq 1 \leq \cdots \leq n\}. \]
	We obtain a category $\Delta$ with objects these finite ordered sets $[n]$ for $n \geq 0$, together with order-preserving functions between them. To get an idea of the structure of this category, let us look at some examples.
	
	There are two possible maps $[0] \rightarrow [1]$, given by $0 \mapsto 0$ and $0 \mapsto 1$.  There is only one map $[1] \rightarrow [0]$, since $[0]$ only has one element.  If we consider maps $[2] \rightarrow [1]$, there are several possibilities.  We could take $0,1 \mapsto 0$ and $2 \mapsto 1$, or $0 \mapsto 0$ and $1,2 \mapsto 1$.  
	
	We could also send all elements of $[2]$ to a single element of $[1]$, but such a map can also be obtained by composing with maps between $[0]$ and $[1]$, so it is not as important to take into account.  Not including such constant maps, we get three maps $[1] \rightarrow [2]$: $0 \mapsto 0$ and $1 \mapsto 1$, $0 \mapsto 1$ and $1 \mapsto 2$, and $0 \mapsto 0$ and $1 \mapsto 2$.
	
	In general, there are $n+1$ generating maps $[n-1] \rightarrow [n]$ and $n+1$ generating maps $[n+1] \rightarrow [n]$.  So, we can visualize the category $\Delta$ as
	\[ \xymatrix@1{[0] \ar@<.5ex>[r] \ar[r] & [1] \ar@<.5ex>[l] \ar[r] \ar@<.5ex>[r] \ar@<1ex>[r] & [2] \ar@<.5ex>[l] \ar@<1ex>[l] \ar[r] \ar@<.5ex>[r] \ar@<1ex>[r] \ar@<1.5ex>[r] & \cdots. \ar@<.5ex>[l] \ar@<1ex>[l] \ar@<1.5ex>[l]} \]
\end{example}

Observe that the category $\Delta$ that we have just described closely resembles the ``shape" of a simplicial set, as defined in the previous section, except that the morphisms are going in the wrong direction.  

From a categorical point of view, we can apply the following formal procedure.

\begin{definition}
	Given any category $\mathcal C$, we can define its \emph{opposite category} $\mathcal C^{\op}$ to be the category with the same objects as $\mathcal C$, but for which $\Hom_{\mathcal C^{\op}}(x,y) = \Hom_\mathcal C(y,x)$.  
\end{definition}

In other words, the direction of the arrows is reversed when passing to the opposite category.  We encourage the reader unfamiliar with this definition to check that composition, associativity, and identity morphisms work as desired in the opposite category.

We can thus consider the category $\Deltaop$ opposite to the category of finite ordered sets, which we can depict as
\[ \xymatrix@1{[0] \ar@<-.5ex>[r] & [1] \ar[l] \ar@<-.5ex>[l] \ar@<-.5ex>[r] \ar@<-1ex>[r] & [2] \ar[l] \ar@<-.5ex>[l] \ar@<-1ex>[l] \ar@<-.5ex>[r] \ar@<-1ex>[r] \ar@<-1.5ex>[r] & \cdots. \ar[l] \ar@<-.5ex>[l] \ar@<-1ex>[l] \ar@<-1.5ex>[l]} \]
Thus, a simplicial set as described in Definition \ref{simpsetdef} is a collection of sets $X_n$, corresponding to the objects $[n]$, and functions between them as indicated by the category $\Deltaop$.  

To be more precise, we introduce the notion of a functor.

\begin{definition}
	Let $\mathcal C$ and $\mathcal D$ be categories.  A \emph{functor} $F \colon \mathcal C \rightarrow \mathcal D$ is given by assigning to every object $x$ of $\mathcal C$ an object $Fx$ of $\mathcal D$, and assigning to every morphism $f \colon x \rightarrow y$ in $\mathcal C$ a morphism $F(f) \colon Fx \rightarrow Fy$.  This assignment should respect composition and identities, so $F(g \circ f) = F(g) \circ F(f)$ and $F(\id_x) = \id_{Fx}$.
\end{definition}

We already saw one example of a functor in our discussion of simplicial complexes above.

\begin{example}
	Geometric realization defines a functor $|-| \colon \mathcal{ASC} \rightarrow \Top$.
\end{example}

Now, let us return to to this idea of a simplicial set as some kind of diagram in the shape of $\Deltaop$.  When we have a functor from a small category $\mathcal C$ to a (possibly large) category $\mathcal D$, we often regard it as a diagram in which we take objects and morphisms in $\mathcal D$ as indicated by the category $\mathcal C$.  For example, if $\mathcal C$ is the category
\[ \bullet \rightarrow \bullet \leftarrow \bullet, \]
then a functor from $\mathcal C$ to an arbitrary category $\mathcal C$ is given by a diagram 
\[ d_1 \rightarrow d_2 \leftarrow d_3 \]
in $\mathcal D$.  Specifically, each $d_i$ is an object of $\mathcal D$, and the arrows are morphisms in $\mathcal D$.  Applying this idea to the small category $\Deltaop$ and the category $\Sets$, we can reformulate Definition \ref{simpsetdef} in the following way.

\begin{definition}
	A \emph{simplicial set} is a functor $K \colon \Deltaop \rightarrow \Sets$.
\end{definition}

In particular, we have packaged all the information about face and degeneracy maps into the structure of the category $\Deltaop$.  While the first definition we gave is perhaps more accessible initially, this alternative approach gives an efficient way to work with simplicial sets.

To illustrate the utility of this definition, let us give a new way to define the $n$-simplex $\Delta[n]$.

\begin{definition}
	The $n$-\emph{simplex} $\Delta[n]$ is the simplicial set given by
	\[ \Delta[n]_k = \Hom_\Delta([k], [n]), \]
	where the right-hand side denotes the set of functors from the category $[k]$ to the category $[n]$.  The face maps $\Delta[n]_k \rightarrow \Delta[n]_{k-1}$ are induced from the corresponding map $[k-1] \rightarrow [k]$ in the category $\Delta$, and similarly for the degeneracy maps.  We can write this definition more efficiently as
	\[ \Delta[n] = \Hom_\Delta(-, [n]). \]
\end{definition}

This definition is an example of what is called a \emph{representable functor} in category theory, namely a functor defined by taking morphisms into a fixed object.  Many nice properties of $\Delta[n]$ can be deduced from the fact that it is a representable functor.

Just as we have done for several other mathematical structures, we would like to think of a category of simplicial sets, for which we need an appropriate notion of morphism between simplicial sets.  Intuitively, such a morphism $K \rightarrow L$ should be given by a sequence of set functions $K_n \rightarrow L_n$ for $n \geq 0$, commuting appropriately with all face and degeneracy maps.  Using the categorical definition of a simplicial set as a functor, we can get a more concise definition.  To do so, we introduce the following general definition from category theory.

\begin{definition} \label{nattrans}
	Let $\mathcal C$ and $\mathcal D$ be categories and $F, G \colon \mathcal C \rightarrow \mathcal D$ functors between them.  A \emph{natural transformation} $\eta \colon F \Rightarrow G$ is given by, for each object $x$ of $\mathcal C$, a function $\eta_x \colon F(x) \rightarrow G(x)$, such that, for any function $f \colon x \rightarrow y$ in $\mathcal C$ the diagram
	\[ \xymatrix{F(x) \ar[r]^{\eta_x} \ar[d]_{F(f)} & G(x) \ar[d]^{G(f)} \\
		F(y) \ar[r]^{\eta_y} & G(y)} \]
	commutes.
\end{definition}

Applying this definition to simplicial sets, we obtain the following.

\begin{definition}
	A morphism $K \rightarrow L$ of simplicial sets is given by a natural transformation from $K$ to $L$.
\end{definition}

In other words, the data of such a morphism is given by set maps $K_n \rightarrow L_n$ for each $n \geq 0$ that commute appropriately with the face and degeneracy maps, as specified by the diagrams in Definition \ref{nattrans}.

We denote the category of simplicial sets, with natural transformations as morphisms, by $\SSets$.  Observe that there is a natural notion of isomorphism of simplicial sets, where each component $K_n \rightarrow L_n$ is an isomorphism of sets.

Thinking about a simplicial set as a functor, rather than in terms of the more geometric descriptions we have given previously, can take some getting used to, and the description of morphisms can be still more mysterious.  While the translation between the two is formalized by geometric realization, which we describe in the next section, the key to working with simplicial sets is in keeping the two ways of thinking about them in mind simultaneously.

\section{The relationship between simplicial sets and topological spaces} \label{ssetstop}

Our initial motivation for defining simplicial sets was to provide a good combinatorial model for topological spaces.  Now that we have some categorical language available, we can say more about the relationship between the categories $\Top$ and $\SSets$.

In the previous section, we claimed that the geometric realization of an abstract simplicial complex is an example of a functor.  We now describe how to upgrade this construction to a functor from the category of simplicial sets to the category of topological spaces.

\begin{definition}
	Let $K$ be a simplicial set.  Its \emph{geometric realization} $|K|$ is defined by
	\[ |K| = \coprod_{n \geq 0} K_n \times \Delta^n / \sim, \]
	where the identification is given by gluing the topological simplices $\Delta^n$ as specified by the face maps of $K$, and identifying degenerate simplices to the nondegenerate simplices from which they arise.
\end{definition}

Observe that this definition is somewhat more complicated than the definition that we gave for the geometric realization of an abstract simplicial complex, but it follows the same heuristic idea.  We invite the reader to revisit previous examples and verify that this definition gives the expected topological output.

We can also describe an explicit functor taking a topological space to a simplicial set.  

\begin{definition}
	Let $X$ be a topological space.  Its \emph{singular simplicial set} $S(X)$ is defined by $S(X)_n = \Hom_{\Top}(\Delta^n, X)$. 
\end{definition}

Readers familiar with singular homology might recognize $S(X)_n$ as the generating set for the free abelian group of singular $n$-chains on $X$.  We leave it as an exercise to verify that the singular simplicial set construction is functorial.

Observe that we have deviated from our treatment of simplicial complexes here, as the singular simplicial set is a much more complicated object than the ``abstractification" of a geometrically-described simplicial complex.  However, this construction is much more general, as it can be applied to any topological space, not just one that has been triangulated so as to have the structure of a simplicial complex.

We would like to understand the result of applying these two functors consecutively, for which we review the definition of adjoint functors.

\begin{definition}
	Let $\mathcal C$ and $\mathcal D$ be categories and $F \colon \mathcal C \rightarrow \mathcal D$ and $G \colon \mathcal D \rightarrow \mathcal C$ functors between them.  The pair $(F,G)$ is an \emph{adjoint pair} of functors if, for any object $X$ of $\mathcal C$ and object $Y$ of $\mathcal D$, there is a natural bijection of sets
	\[ \Hom_\mathcal D(FX,Y) \cong \Hom_\mathcal C(X,GY). \]
\end{definition}

Adjoint functors appear frequently throughout mathematics, even when not always labeled as such.  One example is the tensor-hom adjunction for modules; we refer the reader to \cite[\S 4.1]{riehl} for a detailed list of many others.

\begin{prop} \cite[I.2.1]{gj}
	The singular functor $\SSets \rightarrow \Top$ is right adjoint to geometric realization.
\end{prop}

We can say even more about the output of the singular functor, for which we need the following definition.  Recall the horns $\Lambda^k[n]$ from Definition \ref{ssetexs}.

\begin{definition} \label{kancxdefn}
	A \emph{Kan complex} is a simplicial set $K$ such that for any $n \geq 1$, $0 \leq k \leq n$, and any diagram of the form
	\[ \xymatrix{\Lambda^k[n] \ar[r] \ar[d] & K \\
		\Delta[n] \ar@{-->}[ur] & } \]
	there is a dotted arrow lift making the diagram commute.
\end{definition}

We claim that Kan complexes are the simplicial sets that most effectively model topological spaces.  In our development of simplicial sets from simplicial complexes, we found we needed to introduce an orientation on simplices to make products behave correctly.  However, this information is lost when we apply geometric realization.  

What happens when we take the singular functor?  For example, given a topological space $X$, 1-simplices in $S(X)$ are given by maps from $\Delta^1$ to $X$, namely, paths in $X$.  But paths $\gamma \colon \Delta^1 \rightarrow X$ in a topological space are only oriented in that we typically parametrize $\Delta^1$ as the interval $[0,1]$; the image of such a map does not have any orientation as a subspace of $X$.  In particular, the path with the reverse orientation also specifies a 1-simplex of $S(X)$ with the same image.  Thus, a 1-simplices of a simplicial set that arises from the singular functor should have ``inverse" 1-simplices.

Recall that the horn $\Lambda^0[2]$ can be depicted as
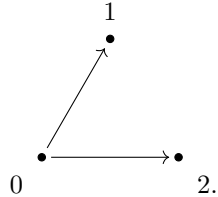
\[ \begin{tikzpicture}[scale=0.7]
	\foreach \tet/\lab/\nom/\pos in {-30/A/2./below right,90/B/1/above,210/C/0/below left}{
		\draw[fill=black] (\tet:1.5) circle (0.07);
		\path (\tet:1.5) node[label=\pos:$\nom$] (\lab) {};
	}
	\draw[->] (C) to (B);
	\draw[->] (C) to (A);
\end{tikzpicture} \]
Consider the map $\Lambda^0[2] \rightarrow S(X)$ that sends the edge $0 \rightarrow 1$ to the 1-simplex corresponding to $\gamma$ and the edge $0 \rightarrow 2$ to the constant path at the initial vertex of the path.  We can fill in this horn by taking the remaining 1-simplex to be the one corresponding to the reverse path of $\gamma$, and taking the 2-simplex given by the map $\Delta^2 \rightarrow X$ to be a homotopy from the concatenation of these two paths to the constant path at the appropriate endpoint.  Observe that this filling is not unique, since there are many choices of such a homotopy, but we do not require uniqueness of the lifts in the definition of Kan complex. 

Indeed, similar arguments can be made in higher dimensions, and we have the following result.

\begin{prop}
	Let $X$ be a topological space.  Then the simplicial set $S(X)$ is a Kan complex. 
\end{prop}

We revisit Kan complexes in the next section, and provide some further intuition there.  In addition to the result we have just stated, further motivation for Kan complexes is indicated in Proposition \ref{kancxwhe} below.

To say more about relationship between the geometric realization and singular functors, we recall the following notions from classical homotopy theory of topological spaces.

\begin{definition}
	Let $X$ be a topological space and $x_0 \in X$ a specified basepoint.  For any $n \geq 0$, we define 
	\[ \pi_n(X, x_0) = [S^n, X]_*, \]
	where the right-hand side denotes the set of homotopy classes of basepointed maps from the $n$-sphere (itself equipped with a chosen basepoint) to $X$.  When $n=1$, this set has a group structure and is called the \emph{fundamental group}.  When $n \geq 2$, it has the structure of an abelian group.
\end{definition}

In an abuse of notation, we omit the basepoint and simply write $\pi_n(X)$. When $X$ is path-connected there is no ambiguity in doing so; more care needs to be taken in what follows in more general cases, but as the details are more annoying than difficult, we omit them here; see \cite[\S 4.1]{hatcher}.  

In fact, each $\pi_n$ defines a functor, so that we obtain from any continuous map $X \rightarrow Y$ a corresponding $\pi_n(X) \rightarrow \pi(Y)$ that is a function when $n=0$ and a group homomorphism otherwise.

\begin{definition}
	A continuous map $X \rightarrow Y$ of topological spaces is a \emph{weak homotopy equivalence} if the induced maps
	\[ \pi_n(X) \rightarrow \pi_n(Y) \]
	are isomorphisms for all $n \geq 0$.
\end{definition}

This notion, while perhaps less familiar than that of homotopy equivalence, is central to the study of the homotopy theory of topological spaces.  For nice spaces, such as CW complexes, the two notions agree \cite[4.5]{hatcher}.  

With this definition in hand, we return to the comparison of topological spaces and simplicial sets.

\begin{prop}
	For any topological space $X$, the natural map $\varepsilon_X \colon |S(X)| \rightarrow X$ is a weak homotopy equivalence.
\end{prop}

Using an analogous notion of homotopy groups for simplicial sets, such as described in \cite[\S 9]{friedman} or \cite[\S I.7]{gj}, we also have the following.

\begin{prop} \label{kancxwhe} \cite[I.11.1]{gj}
	For any Kan complex $K$, the natural map $K \rightarrow S|K|$ is a weak homotopy equivalence.
\end{prop}

Roughly speaking, the previous two propositions combine to give the heart of the following result.

\begin{theorem} \cite[I.11.4]{gj}
	The adjoint pair 
	\[ |-| \colon \SSets \rightleftarrows \Top \colon S \]
	defines an equivalence of homotopy theories.
\end{theorem}

We do not go into the details of how to make this theorem, or even the definition of ``homotopy theory", rigorous here; it can be realized, for example, as a Quillen equivalence of model categories, but we do not wish to go into the details of that theory here.  The interested reader can consult \cite{ds}.

One might wonder, after all this work, whether developing the theory of simplicial sets, and proving its equivalence with topological spaces is worth it.  Let us give an example to give a taste of why working with simplicial sets is often much easier than working with topological spaces.

\begin{example} \label{mapping}
	Let $X$ and $Y$ be topological spaces.  One can consider $\Map(X,Y)$, the set of all continuous functions $X \rightarrow Y$, as a topological space; the standard method is to use the compact-open topology \cite[\S 46]{munkres}, which does what we want it to do, but is not entirely intuitive.
	
	On the other hand, let $K$ and $L$ be simplicial sets.  We can define a simplicial set $\Map(K,L)$ by $\Map(K,L)_n = \Hom(K \times \Delta[n], L)$, where $\Hom$ denotes the usual set of morphisms in the category of simplicial sets.  This definition is as straightforward as we could possibly wish.
	
	Indeed, this example itself leads to a wealth of methods for working with mapping spaces.  It is much more common to work with ``mapping simplicial sets" than with ``mapping topological spaces" in more general categorical settings.  The notion of a simplicial category, or category equipped with a natural notion of mapping simplicial sets, has been an important one with a wide range of applications.
\end{example}

\section{Simplicial sets and categories} \label{ssetscat}

We have discussed why simplicial sets are important in topology, but they also play an important role in category theory.  Let us begin by looking at a construction that links simplicial sets and categories.  More explicitly, we want to look at a functor from the category $\Cat$, whose objects are small categories and whose morphisms are functors between them, and the category $\SSets$ of simplicial sets, as defined in the previous section.

We can think of the finite ordered set $[n] = \{0 \leq 1 \leq \cdots \leq n\}$ as a category, in which we replace the ordering by arrows, so we have a single arrow $i \rightarrow j$ whenever $i \leq j$.  Indeed, any poset can be regarded as a category in an analogous way.

Given any other (small) category $\mathcal C$, a functor $[n] \rightarrow \mathcal C$ is given by a string of $n$ composable arrows in $\mathcal C$:
\[ x_0 \rightarrow x_1 \rightarrow \cdots \rightarrow x_n. \]

\begin{definition}
	Let $\mathcal C$ be a small category.  The \emph{nerve} of $\mathcal C$ is the simplicial set $\nerve(\mathcal C)$, given in degree $n$ by
	\[ \nerve(\mathcal C)_n = \Hom_{\Cat}([n], \mathcal C). \]
\end{definition}

Note that the 0-simplices of the nerve correspond to the objects of $\mathcal C$, the 1-simplices correspond to the morphisms of $\mathcal C$, and more general $n$-simplices correspond to chains of $n$ composable morphisms in $\mathcal C$.

\begin{example}
	Let $G$ be a group, regarded as a category with a single object.  Then the geometric realization of the nerve, $|\nerve(G)|$ gives a model for the classifying space of $G$; see \cite[1B.7]{hatcher}.
\end{example}

In fact, the nerve construction defines a functor $\Cat \rightarrow \SSets$.  One of the especially nice features of this functor is the following result.

\begin{prop} \label{nerveff}
	The nerve functor $\Cat \rightarrow \SSets$ is fully faithful; in other words, for any categories $\mathcal C$ and $\mathcal D$, there is an isomorphism 
	\[ \Hom_{\Cat}(\mathcal C, \mathcal D) \cong \Hom_{\SSets}(\nerve(\mathcal C), \nerve(\mathcal D)). \]
\end{prop}

In other words, for many purposes we can just as easily work with the nerves of categories as with the categories themselves, since the functions between them are the same.

How can we detect which simplicial sets are nerves of categories?  The key point to recognize is that a nerve of a category has to have ``composition" of 1-simplices.  So, for example, any time we have a configuration 
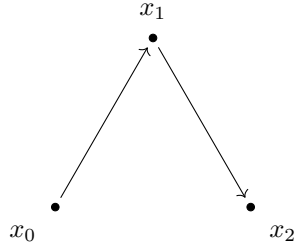
\[ \begin{tikzpicture}
	\foreach \tet/\lab/\nom/\pos in {-30/A/x_2/below right,90/B/x_1/above,210/C/x_0/below left}{
		\draw[fill=black] (\tet:1.5) circle (0.05);
		\path (\tet:1.5) node[label=\pos:$\nom$] (\lab) {};
	}
	\draw[->] (C) to (B);
	\draw[->] (B) to (A);
\end{tikzpicture} \] 
of 1-simplices, there must be another 1-simplex $x_0 \rightarrow x_2$, and these three 1-simplices form the boundary of a 2-simplex.  

Recall that this configuration was described previously as the horn $\Lambda^1[2]$.  We can encode iterated composition, as well as associativity, using higher-dimensional horns.

\begin{prop}
	A simplicial set $K$ is isomorphic to the nerve of a category if and only if for any $n \geq 1$, $0 < k < n$, and any diagram of the form
	\[ \xymatrix{\Lambda^k[n] \ar[r] \ar[d] & K \\
		\Delta[n] \ar@{-->}[ur] & } \]
	there is a unique dotted arrow lift making the diagram commute.
\end{prop}

The horns used here are called \emph{inner horns}, since $k$ does not take the values 0 and $n$.  One might ask what lifts with respect to the \emph{outer horns} give us.  Let us think about the two outer horns when $n=2$:
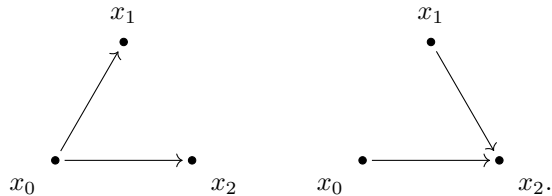
\[ 
\begin{tabular}{ c c }
	\begin{tikzpicture}[scale=0.7]
		\foreach \tet/\lab/\nom/\pos in {-30/A/x_2/below right,90/B/x_1/above,210/C/x_0/below left}{
			\draw[fill=black] (\tet:1.5) circle (0.07);
			\path (\tet:1.5) node[label=\pos:$\nom$] (\lab) {};
		}
		\draw[->] (C) to (B);
		\draw[->] (C) to (A);
	\end{tikzpicture}
	&
	$\quad$
	\begin{tikzpicture}[scale=0.7]
		\foreach \tet/\lab/\nom/\pos in {-30/A/x_2./below right,90/B/x_1/above,210/C/x_0/below left}{
			\draw[fill=black] (\tet:1.5) circle (0.07);
			\path (\tet:1.5) node[label=\pos:$\nom$] (\lab) {};
		}
		\draw[->] (B) to (A);
		\draw[->] (C) to (A);
	\end{tikzpicture}
\end{tabular} \]

Observe that having unique lifts with respect to these horns corresponds to the existence of inverse morphisms, extending the arguments we made in our discussion of Kan complexes in the previous section.

\begin{definition} 
	A category $\mathcal C$ is a \emph{groupoid} if all its morphisms have inverses.
\end{definition}

\begin{prop}
	A simplicial set $K$ is isomorphic to the nerve of a groupoid if and only if for any $n \geq 1$, $0 \leq k \leq n$, and any diagram of the form
	\[ \xymatrix{\Lambda^k[n] \ar[r] \ar[d] & K \\
		\Delta[n] \ar@{-->}[ur] & } \]
	there is a unique dotted arrow lift making the diagram commute.
\end{prop}

If we drop the uniqueness assumption and only ask for the existence of such a lift, we recover the definition of Kan complex from Definition \ref{kancxdefn}.  The approach here suggests that Kan complexes are something like ``groupoids up to homotopy", whereas the results of the previous section suggest that they are the simplicial sets that best model topological spaces.  However, these two perspectives shed light on one another.  For example, in a groupoid the ``directionality" of morphisms is less important than in an ordinary category, since every morphism has an inverse going in the opposite direction.  But topological spaces also lack such directionality, as we saw in our discussion of why the output of the singular functor is a Kan complex.   

It is not unreasonable to apply the same weakening to the simplicial sets that appear as nerves of categories, namely, to drop the uniqueness condition.  The following definition was first made by Boardman and Vogt \cite{bv}.

\begin{definition}
	A \emph{weak Kan complex} or \emph{quasi-category} is a simplicial set $K$ if for any $n \geq 2$, $0 < k < n$, and any diagram of the form
	\[ \xymatrix{\Lambda^k[n] \ar[r] \ar[d] & K \\
		\Delta[n] \ar@{-->}[ur] & } \]
	there is a dotted arrow lift making the diagram commute.
\end{definition}

Perhaps surprisingly, such simplicial sets were not nearly so prominent until more recently, but let us give a taste of why they might be important.  We start with the notion of equivalence of categories; recall the definition of a fully faithful functor from Proposition \ref{nerveff}.

\begin{definition}
	A functor $F \colon \mathcal C \rightarrow \mathcal D$ is an \emph{equivalence of categories} if it is fully faithful and if it is also \emph{essentially surjective}, in the sense that any object of $\mathcal D$ is isomorphic to an object in the image of $F$.
\end{definition}

\begin{prop} \cite[3.1.3]{book}
	Let $\mathcal C \rightarrow \mathcal D$ be an equivalence of categories.  Then the induced functor on nerves 
	\[ \nerve(\mathcal C) \rightarrow \nerve(\mathcal D) \]
	is a weak homotopy equivalence.
\end{prop}

However, the converse of this proposition is false.  

\begin{example}
	Consider the category $[0]$, with one object and only an identity morphism, and the category $[1]$, with two objects and a single morphism from one to the other.  Consider either of the two functors
	\[ [0] = (\bullet) \rightarrow (\bullet \rightarrow \bullet) = [1]. \]
	Neither functor is an equivalence because the object not in the image is not isomorphic to the object that is.
	
	However, the nerve of $[0]$ is the 0-simplex $\Delta[0]$, and the nerve of $[1]$ is the 1-simplex $\Delta[1]$.  Since the geometric realization of each of these simplicial sets results in a contractible space, the above functor induces a weak homotopy equivalence on nerves.
\end{example}	

The phenomenon that we are observing in this example is that weak homotopy equivalences are determined by what happens after applying geometric realization.  In particular, to return to a theme introduced above, we lose all directionality of simplices, and more specifically, we do not retain information about whether two objects in the original category were isomorphic to one another or simply just had a morphism between them.

We thus return to the idea what topological spaces are more accurately modeled by Kan complexes, which behave more like groupoids, rather than by arbitrary simplicial sets.  If we put quasi-categories, rather than Kan complexes, into the main role, then we get a related notion of ``weak equivalence" that is more reflective of equivalences of categories rather than of weak homotopy equivalences of topological spaces.  The theory here is quite deep, so we do not go into the details here, but offer it as a motivation for the study of quasi-categories.

\section{Beyond simplicial sets} \label{simpobj}

In Example \ref{mapping}, we illustrated one advantage of working with simplicial sets rather than topological spaces.  In this final section of the paper, we describe some further further motivation for the methods used in this paper.  We begin with the notion of more general simplicial objects.

\begin{definition}
	Let $\mathcal C$ be a category.  A \emph{simplicial object in} $\mathcal C$ is a functor $\Deltaop \rightarrow \mathcal C$.
\end{definition}

Let us look at some examples.

\begin{example}
	A \emph{simplicial group} is a simplicial object in $\Gp$, the category of groups and group homomorphisms.  Thus, it consists of groups $G_n$ for each $n \geq 0$ and face and degeneracy maps between them that are required to be group homomorphisms.
	
	In analogy to the comparison between simplicial sets and topological spaces, simplicial groups model \emph{topological groups}, which are topological spaces equipped with the structure of a group in a compatible way.  Describing that compatibility is not too complicated, but not as straightforward as taking a functor $\Deltaop \rightarrow \Gp$.
\end{example}

\begin{example}
	We can similarly take simplicial objects in $\Ab$, the category of abelian groups and group homomorphisms.  While we can think of such  simplicial abelian groups as combinatorial models for topological abelian groups, as in the previous example, we can say more in this case.  An important classical result in the study of simplicial sets is the Dold-Kan Theorem, which states that the category of simplicial abelian groups is equivalent to the category of non-negatively graded chain complexes of abelian groups \cite[III.2.3]{gj}.  Roughly speaking, this equivalence takes a simplicial abelian group to a chain complex whose differentials are given by an alternating sum of the face maps; the reverse map is more complicated.  This result brings simplicial objects into close proximity with homological algebra.
\end{example}

We can similarly model other topological spaces equipped with algebraic structure by taking simplicial objects in the category of the appropriate algebraic objects, for example topological rings by simplicial objects in the category of rings.

From a different point of view, simplicial sets often form a good target category for presheaves, which are simply categories of functors; one can require additional compatibility conditions to get sheaves. Ordinary presheaves are just functors from category to the category $\Sets$ of sets, but they can also be taken simplicially.

\begin{definition}
	Let $\mathcal C$ be a small category.  A \emph{simplicial presheaf on} $\mathcal C$ is a functor $\mathcal C^{\op} \rightarrow \SSets$.  
\end{definition}

While one could just as easily take functors into the category $\Top$, one advantage of the category $\SSets$ here is the fact that its objects are themselves functors.  Thus, a simplicial presheaf $\mathcal C^{\op} \rightarrow \SSets$ can equally be regarded as a functor $\mathcal C^{\op} \times \Deltaop \rightarrow \Sets$.  

While we do not go into the details here, the theory of simplicial presheaves has played an important role in motivic homotopy theory; for a survey see \cite{levine}.

We can bring these two different generalizations together in the study of \emph{bisimplicial sets}, or functors $\Deltaop \rightarrow \SSets$.  We can think of them either as simplicial presheaves on $\Delta$, or simplicial objects in the category of simplicial sets.  Unlike simplicial objects in categories of algebraic objects, there is no immediate analogy to topological spaces with some kind of additional structure.  While we have discussed quasi-categories in more detail above, several different approaches to ``categories up to homotopy" or, more formally, $(\infty,1)$-categories, are given by bisimplicial sets with additional structure.  We refer the reader to \cite{survey} for more details.

\end{document}